# The Markov model of a dynamic system based on experimental data for control problems of bionic prostheses

Meshchikhin I.A., Minkov S., Lichkunov A.A.

**Abstract.**
This article is devoted to a prediction of gait kinematics based on the Markov chains in the delay space. We use hypercubic grid in the delay space of angles in the sagittal plane in the knee and hip joints to construct Markov states. Our experimental signal (obtained from the telemetry of the prosthesis) provides transition probabilities. The resulting Markov model seems to be helpful for estimation of attractor dimension, characteristic frequencies of the system, and so on.

## Introduction

The prosthetic engineering is highly important as a precedent for the modern-style industry:
1. A prosthesis as a product is inseparable from a rehabilitation service.
2. A prosthesis as a product of collaborative robotics exists in a permanent and without-fulltime-control interaction with a human client.
3. A control system has the ability to adapt itself with help of both a prosthesis telemetry data, as well as statistical information about the gait of other patients (Big Data) and can use correction from a rehabilitation therapist within the framework of telemedicine and "Internet of things" concepts.

Thus, the path of evolution of bionic prostheses is a template for a wide class of industries of the future.

Modern bionic knee prostheses have the ability to actively control the moment of resistance in the knee joint. The problem of building a control system for bionic prostheses lies in the absence of a model of the bionic component of the mechatronic system. The construction of a structural [1] model of human locomotion seems to be too complicated task for practical usage. On the other hand, the phenomenological [2] model should be interpretable for an expert correction. This paper deals with the construction of a Markov model of human locomotion on a knee prosthesis. Among the advantages of the Markov model [3] it is worth highlighting the simplicity of construction, the interpretability of the matrix coefficients and an extensive methodological base for its analysis. It is also worth noting that one can move to continuous time and describe the dynamics of the system in the form of Kolmogorov's differential equations [4].

To construct a Markov model, a procedure is required that matches a set of discrete states to the flow of sensory information. It is reasonable to define as states the neighborhoods of characteristic points, or some finite elements in the phase space of a dynamical system. In this work, the delay space is interpreted in an extended way: in addition to the quantity and its delay shifts, it can include other sensory information from the telemetry of the product, forces, moments, values of the quantity taken with a delay [5].

## Feature points and phase space

To provide a metric on a delay space, one should make its axes dimensionless. Since the coordinates x in this space are data from measuring instruments, it is rational to take into account their metrology by referring the values $\bar{x}$ of the measurement results to the corresponding errors $X_\Delta$

$$x = \frac{\bar{x}}{X_\Delta}$$

For a hypercubic grid, dimensionless coordinates are defined as:
$x = \left[\frac{\bar{x}}{X_\Delta} \frac{L}{(\bar{x}_{max} - \bar{x}_{min})}\right]$, where [ ] is rounding, $L(\bar{x}_{max} - \bar{x}_{min})$ is the number of grid cells, is the signal range .

Then the choice of characteristic points $X_i$, the distance between which is less than 1, does not add additional information about the dynamics of the system.

We propose the following rule for selecting characteristic points:
if $R_i < R_0 \approx 1 \; \forall i \in [0, n], R_i = (X_i - x)(X_i - x)^T$ (where n is the sample size), then the point with coordinates x is included in the list of characteristic ones $X_{n+1}$.

The constructed grid contains a maximum of information about the dynamics of the system under study. On the other hand, with a limitation on the maximum value of the measured value, the amount of information [7] can be estimated as:

$$I = \sum_{i=1}^{n} Ln\left(\frac{Xmax_i}{X\Delta_i}\right).$$

A hypercubic partition allows us to estimate the dimension of a dynamical system as:

$N = \frac{n}{2}$, where n is the number of characteristic points, such that:

$$R_{i,j} = (X_i - X_j)(X_i - X_j)^T < R_0 k, k \approx 1.4$$

Visualization of the relationship between the number of neighbors and dimension is shown in Figure 1

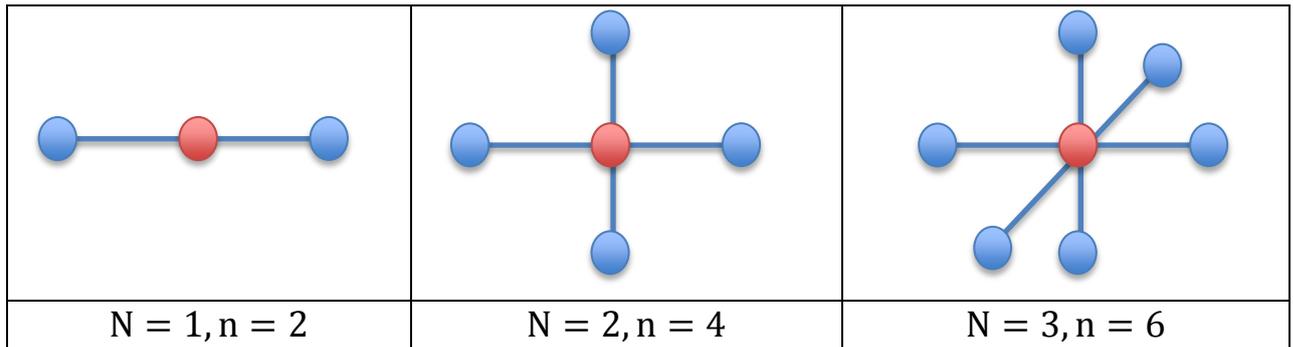

| $N = 1, n = 2$ | $N = 2, n = 4$ | $N = 3, n = 6$ |

Fig. 1 - characteristic points and dimension

## Markov model

The Markov model describes the dependence of the probabilities of belonging to the i-th state at the moment of time t+Δt as:

$$P_{t+\Delta t} = M P_t$$

Elements of the matrix of the Markov model $M_{i,j}$ are composed of frequencies of transition for a period Δt of coordinates *x* from the i-th neighborhood to the j-th. According to the central limit theorem [8], the estimate of the coefficients *M* is corresponding probabilities with an increase in the sample size with the time of training the model. The convergence rate can be estimated by the Chebyshev inequality [9].

## Modal analysis of the Markov model

To be able to analyze and reduce the system, we consider the eigenvalues λ and forms Φ of the matrix $M = \Phi \lambda \Phi^T$ (we are talking about a generic case without nontrivial Jordan cells).

According to the Perron theorem [10], the eigenvalue $\lambda_{max}$ with maximal absolute value of the matrix M is a real positive eigenvalue, and for the particular case of the Markov model, it is equal to one.

The form corresponding to real eigenvalues λ describes dynamics on the attractors [11] of a dynamical system (stationary distributions [12]).

The number of real eigenvalues λ is a lower bound for the degree of the differential equation generating the system.

In generic case the eigenvalues λ of the system are complex.

Complex and negative eigenvalues λ allow us to estimate the oscillation frequency f of the probabilities of being in the i-th neighborhood as:

$$f_i = \frac{1}{n_i \Delta t}, n_i = \frac{2\pi}{\arg(\lambda_i)}.$$

The damping decrement of the mode to the stationary distribution can be defined as:

$$\xi_i = \frac{\ln(|\lambda_{i,i}|)}{2\pi f_i \Delta t}.$$

The modal analysis of the Markov matrix makes it possible to evaluate the dynamics of the system as a superposition of damped oscillation forms.

For the consistency of the model, the condition of the convergence of the model forecast results must be satisfied at $R_0 \to 0$.

With a clear definition of the i-th neighborhood x as the closest to the characteristic point $X_i$, the modal characteristics of the model depend on the choice $R_0$.

With a fuzzy [14] definition of the measure F of belonging to the i-th state as:

$$F_i = \frac{K_i}{\sum K_i}, \text{ where } K_i = \frac{1}{R_i + \alpha} \text{ is the kernel [13].}$$

Then, upon receipt of sensory information, the matrix is M updated with subsequent normalization by the value $F^T F$.

With this filling scheme, the spectrum of the Markov model M is related to the spectrum of the dynamic model by the relationships presented earlier.

## Test case

Let us consider as a test example a linear oscillator with a solution of the form

$$X = Ae^{\xi\omega_d t}\sin(\omega t + \varphi)$$

Let's construct a set of characteristic points in the phase space. The calculation result is shown in Figure 2.

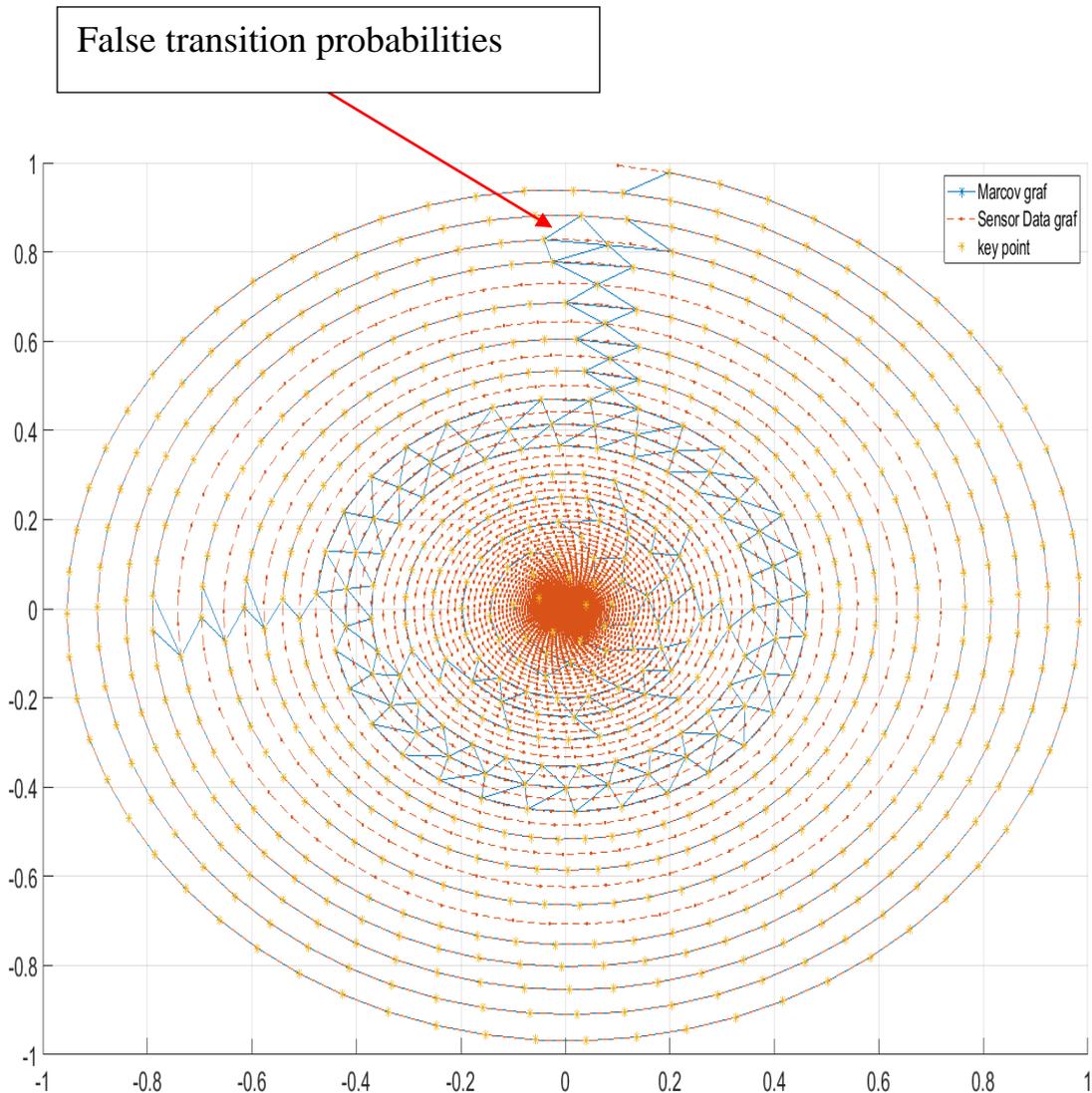

Fig. 2 - characteristic points, the trajectory of damped oscillations and the graph of the Markov chain

The obtained coefficients $M_{i,j}$ are estimates of the probabilities of transition from the i-th state to the j-th state. The probability distribution is shown in Figure 3.

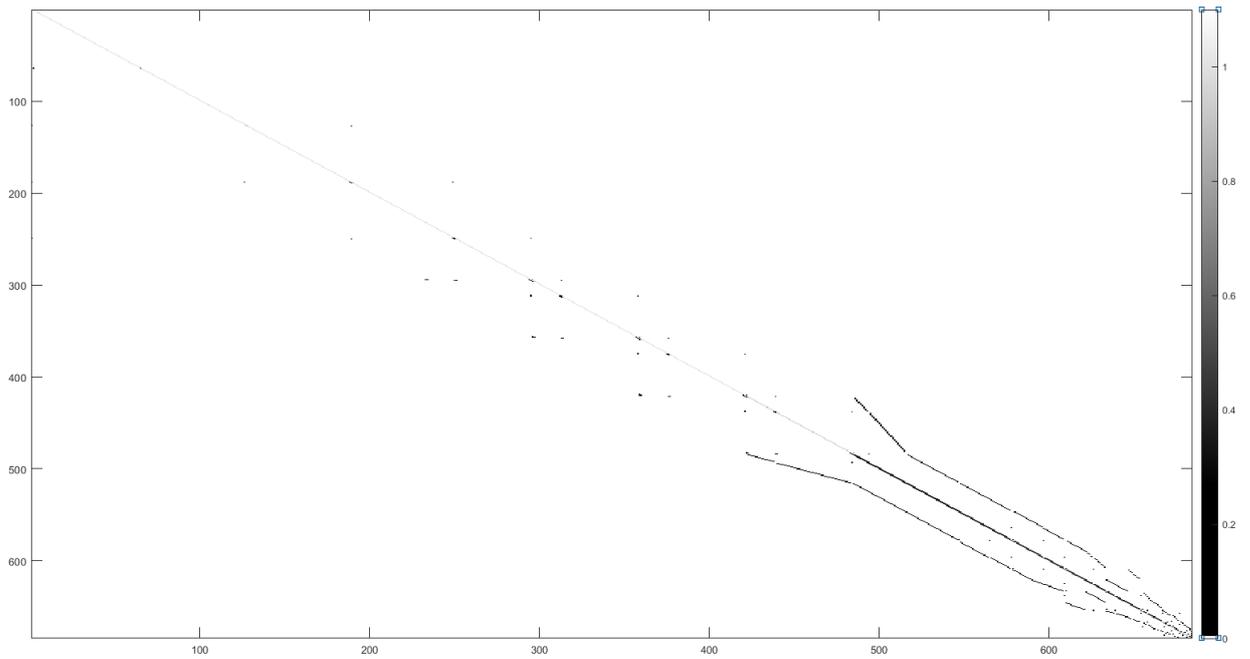

Fig. 3. – Matrix **M**

Let us define the initial conditions as $P_0$. Then the forecast of the state of the system through $n\Delta t$ can be defined as $M^n P_0$. Due to the presence of fictitious transitions due to the division of the phase space into a finite number of elements (and the identification of trajectories passing through them), the resulting forecast will tend to the attractor faster than the real process $P = MP_0$. Nevertheless, if at each iteration all components of the vector are zeroed, except for the maximum N + 1 (N is the dimension of the dynamical system), this effect can be neutralized. Zeroing out more than N + 1 components can result in dummy cyclic paths.

The calculation result is shown in Figure 4 in comparison with the initial data on which the Markov chain was trained.

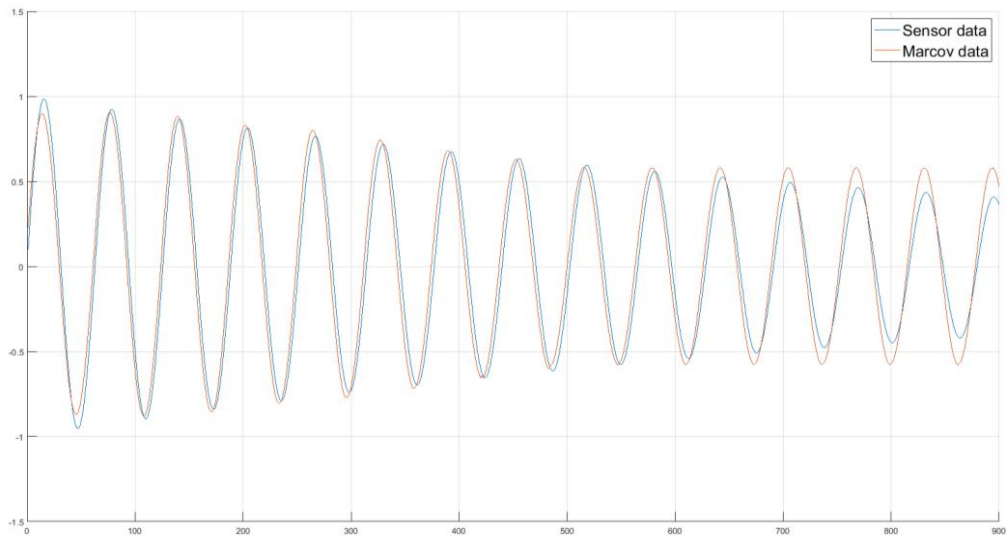

Fig. 4 – The calculation result is presented in comparison with the initial data.

For the model example, the maximum number of "neighboring" characteristic points is 4 – dimension 2.

**Test case modal analysis**

Let us estimate the number of attractors and the degree of nonlinearity of the dynamical system as the number of matrix M eigenvalues λ close to one.
Figure 5 shows the complex eigenvalues λ of the system reconstructed from the data.

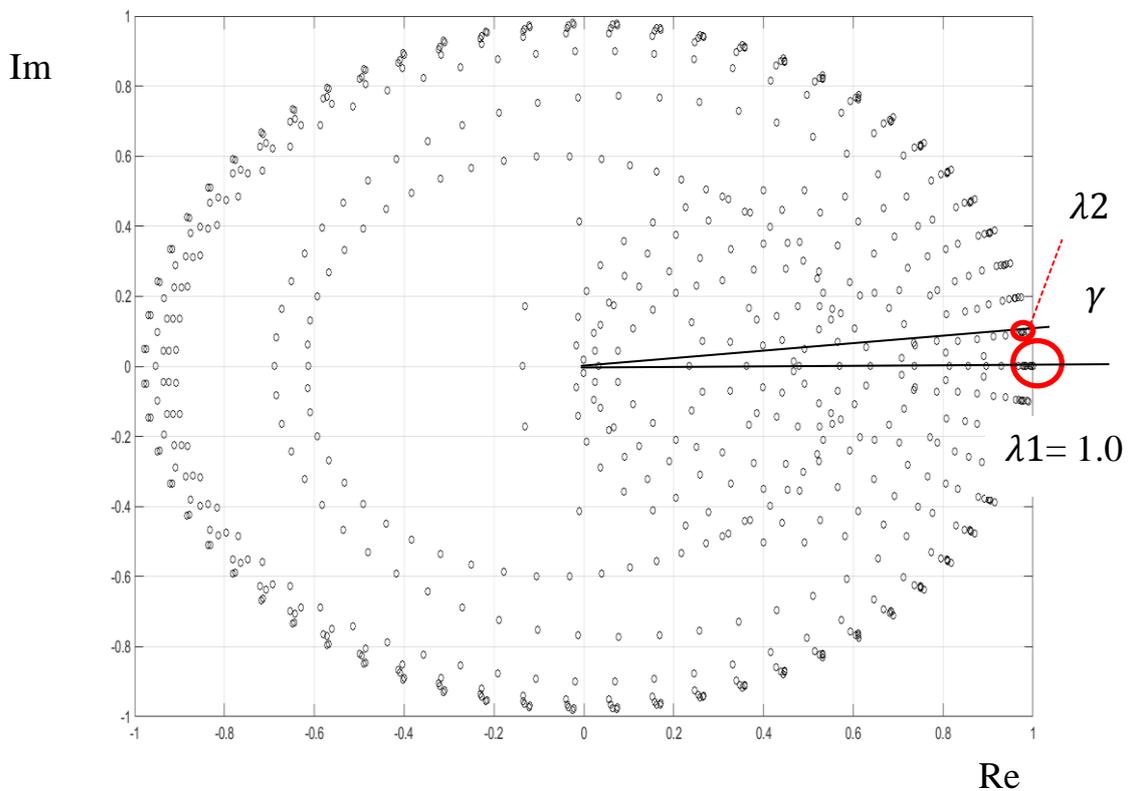

Fig. 5 – eigenvalues λ

Analysis of the eigenvalues λ of the matrix M makes it possible to estimate the characteristic periods of oscillations of a dynamic system:

$$T = \frac{\gamma}{2\pi}\Delta, \text{sec.}$$

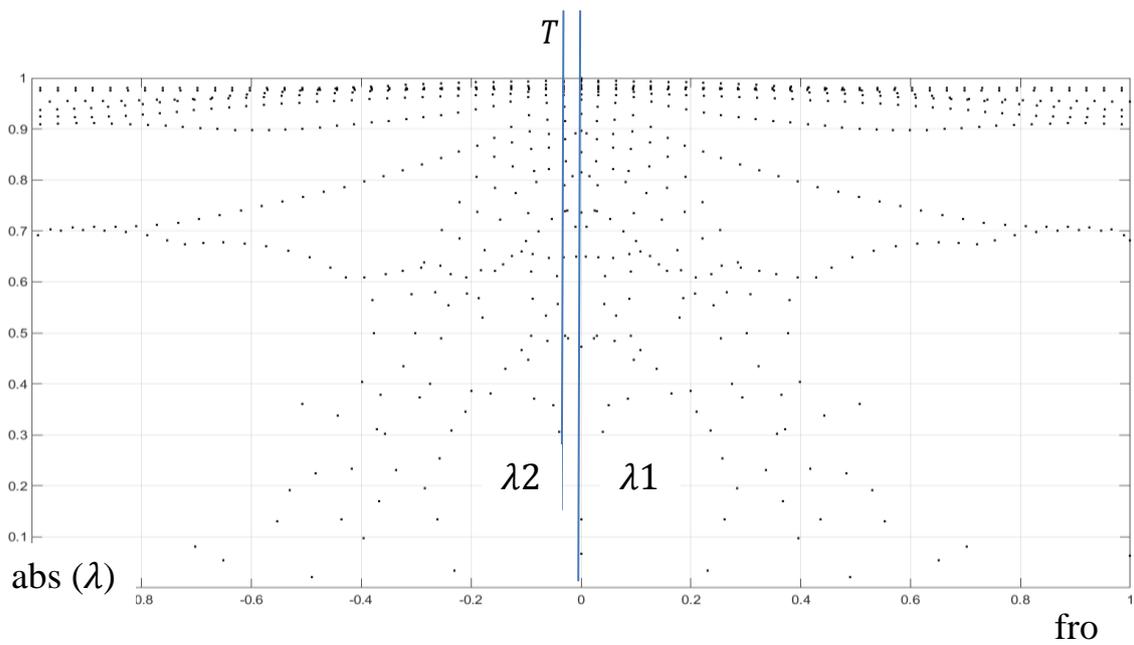

Fig. 6 - Oscillation periods

The closeness of the oscillation periods $kT, k \in N$ to k indicates the linearity of the dynamic system.

The number of attractors of the observed system is 1. Figure 7 shows the eigenvalues form modes for $\lambda_{1,2}$ probabilities in the phase space.

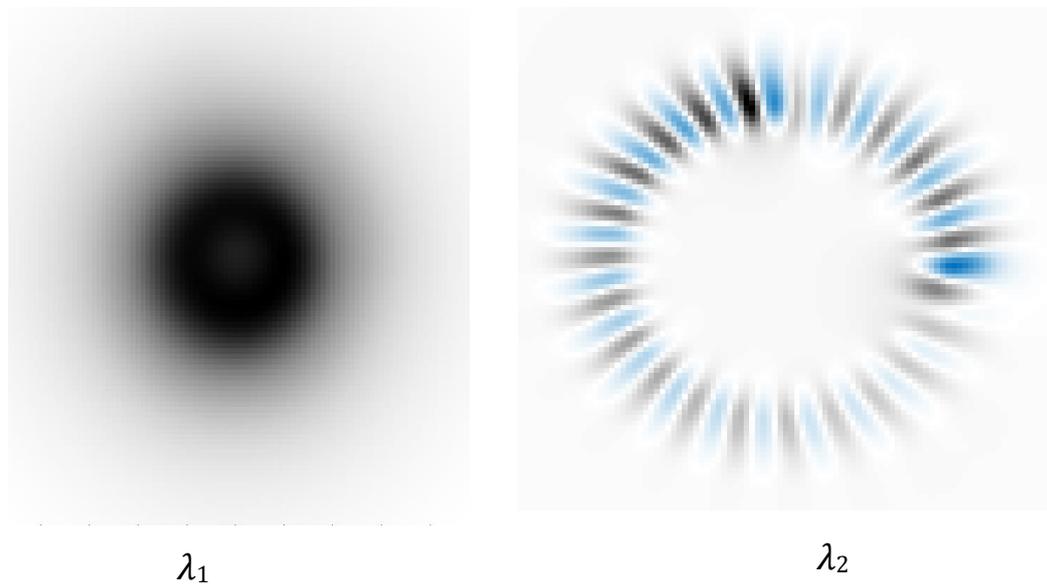

$\lambda_1$ $\lambda_2$

Fig. 7 - eigenforms

The figure 7 shows, that the form corresponding to the eigenvalue equal to one is concentrated in the vicinity of the attractor. The eigenforms corresponding to the eigenvalues $\lambda_{1,2}$ with a nonzero phase reflect the stratification of the system dynamics into a set of oscillatory processes.

## Analysis of telemetry data from the knee prosthesis

Let us represent the dynamics of the walking process in the coordinates of the angles in the knee $\alpha_{KC}$ and the hip $\alpha_{TC}$ joint in Figure 8.
In addition to these angles, the sensory information contains the values $\lambda_{1,2}$ of the angle of the hip joint, taken with a lag. To build a model, consider the trajectory of a person's gait in a 4-dimensional space of angles, as well as angles taken with lags.

A diagram of the measured angles is shown in Figure 8.

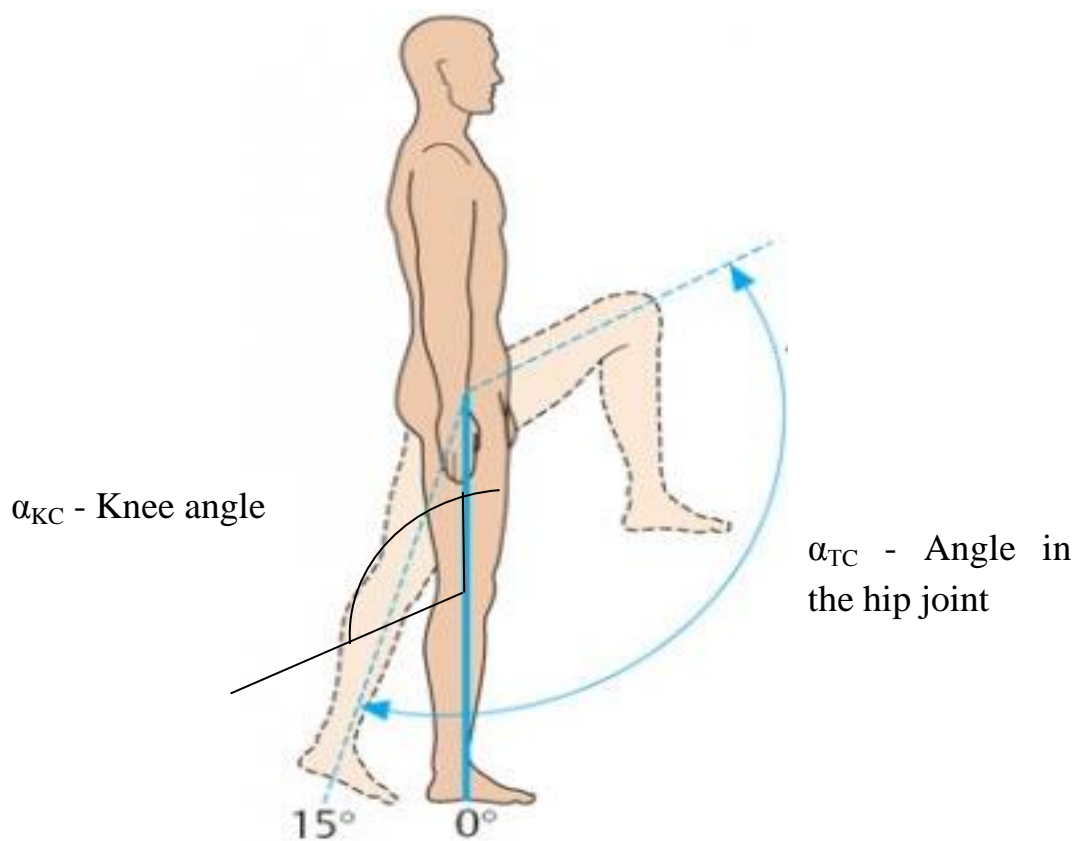

$\alpha_{KC}$ - Knee angle

$\alpha_{TC}$ - Angle in the hip joint

Fig. 8 - measured angles

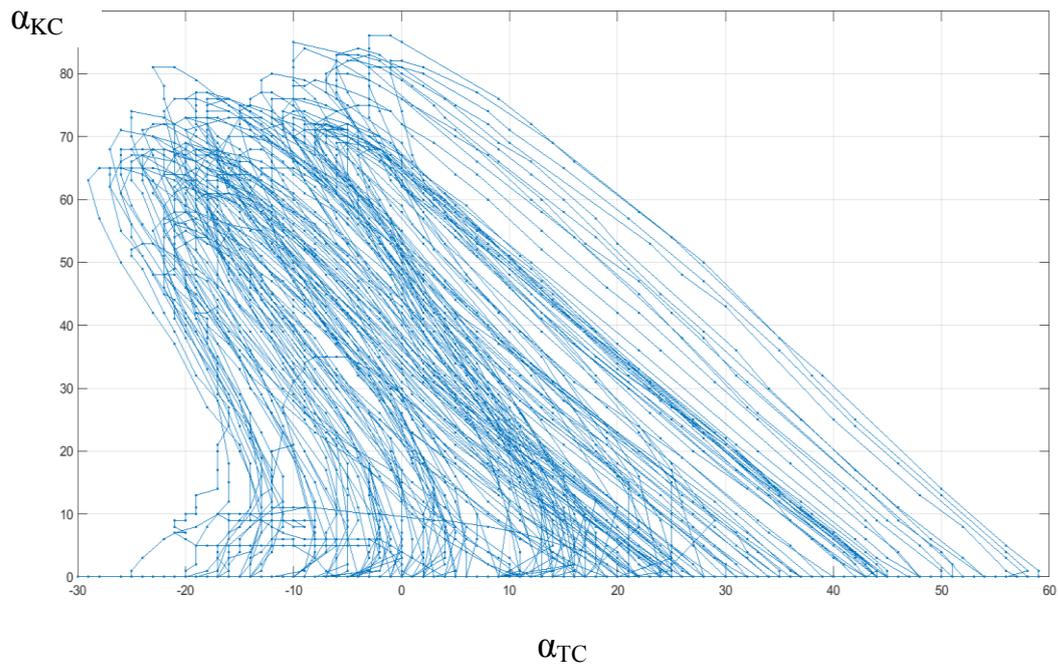

Fig. 9 - Kinematics of gait of an amputated person

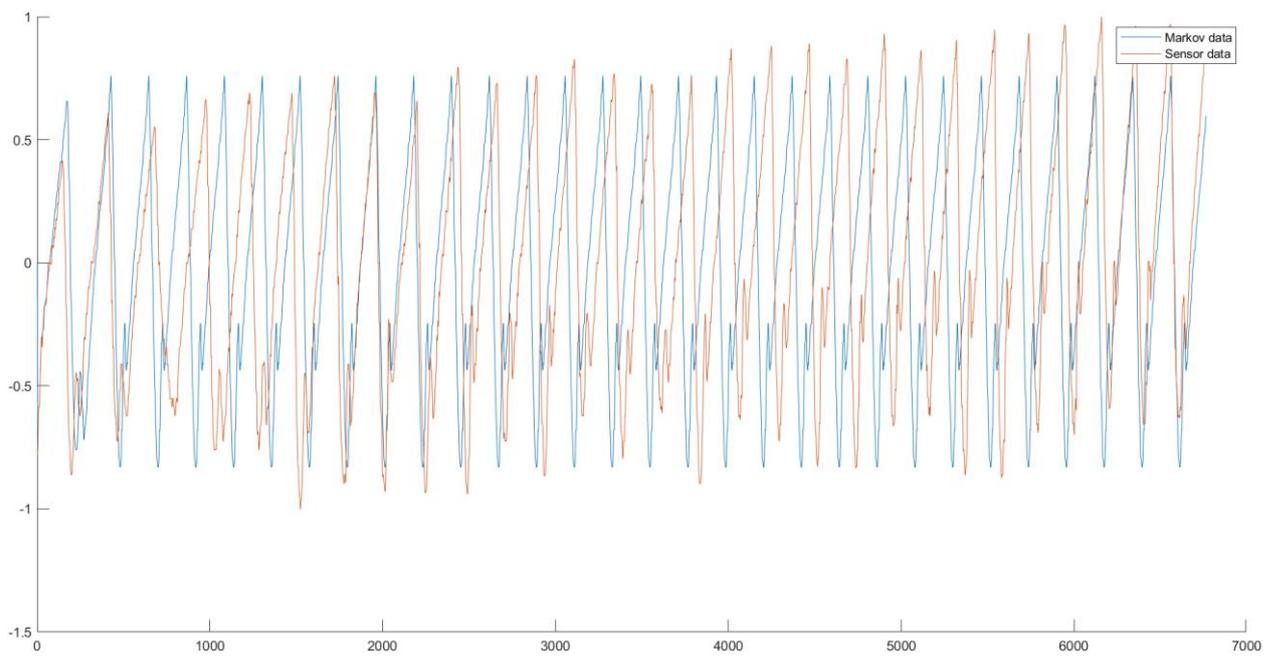

Fig. 10a - Kinematics of gait of an amputated person (normalized angle of the hip joint)

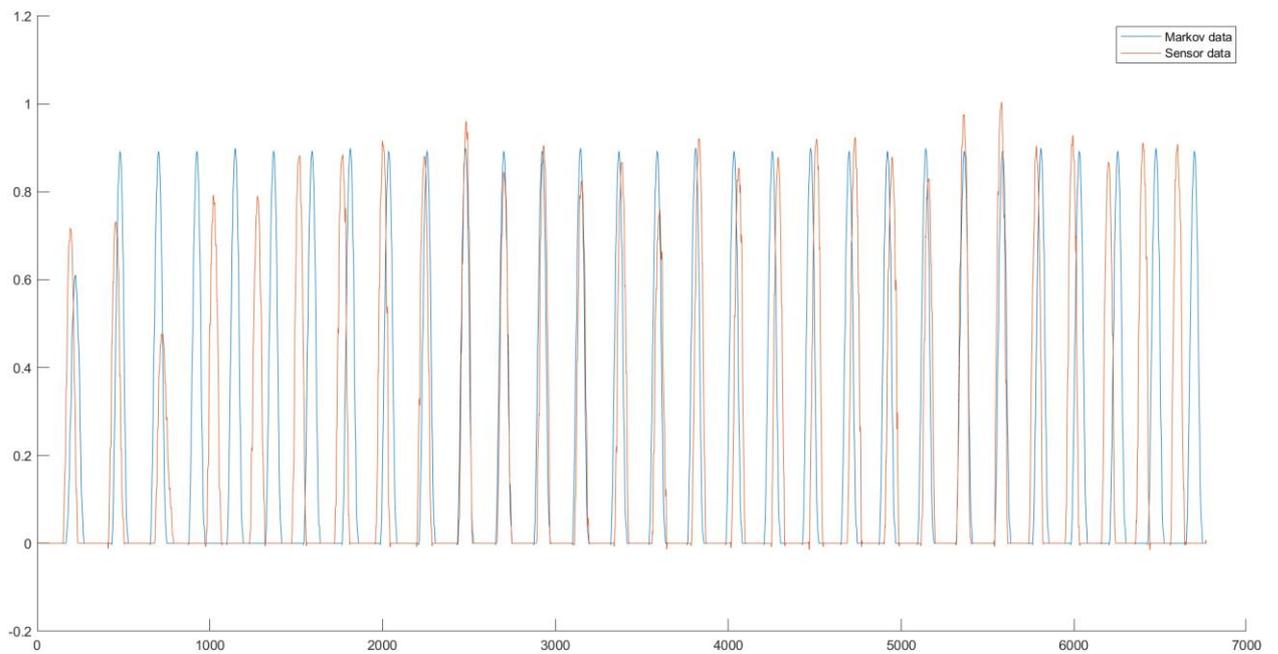

Fig. 10b - Kinematics of gait of a person with amputation (normalized angle of the knee joint)

The dimension was calculated for the case of a uniform hypercubic grid. The dimension estimation graph is shown in Figure 11.

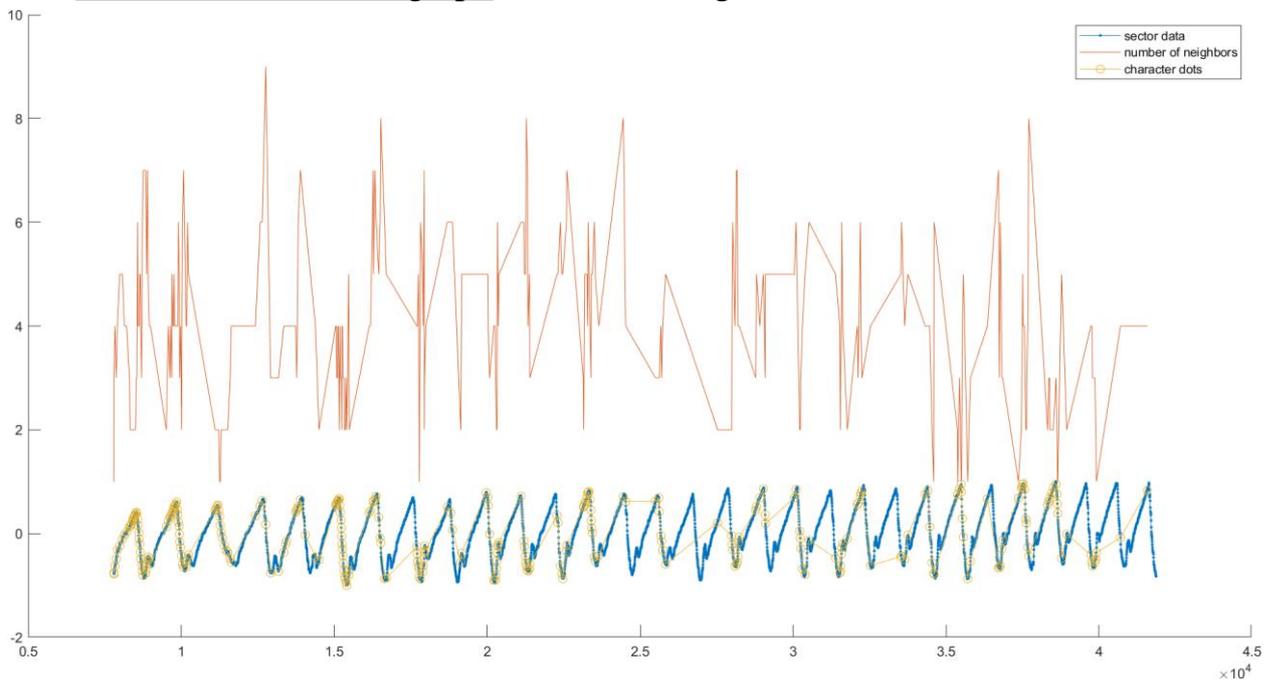

Fig. 11 - Graph of gait kinematics of a person with amputation (normalized angle of the hip joint) with overlapping characteristic points and a graph of the number of neighbors for each characteristic point.

As seen from Fig. 11, the number of neighbors lies in the range [1,8] (a single overshoot towards 9 can be neglected due to discreteness and unevenness of data).

This indicates the correct choice of the dimension of our system. In the case of a continuous signal, all points have a number of neighbors equal to 2 * N = 8. However, due to the discreteness of the signal, as well as the unevenness of the location of the points, we get such a range in the number of neighbors. This moment can be corrected by interpolating additional signal points.

Also here we can talk about the criterion for choosing R0. If most of the points had less than 2 neighbors, this would indicate that there are few characteristic points in the system, which means that the Markov model is under-trained. In this system, 78.4% of points with neighbors are more than 2, which indicates the correct R0.

Modal analysis of prosthetic kinematics

Figure 12 shows the periods of oscillations obtained from the analysis of the eigenvalues $\lambda$ of the matrix of the Markov model.

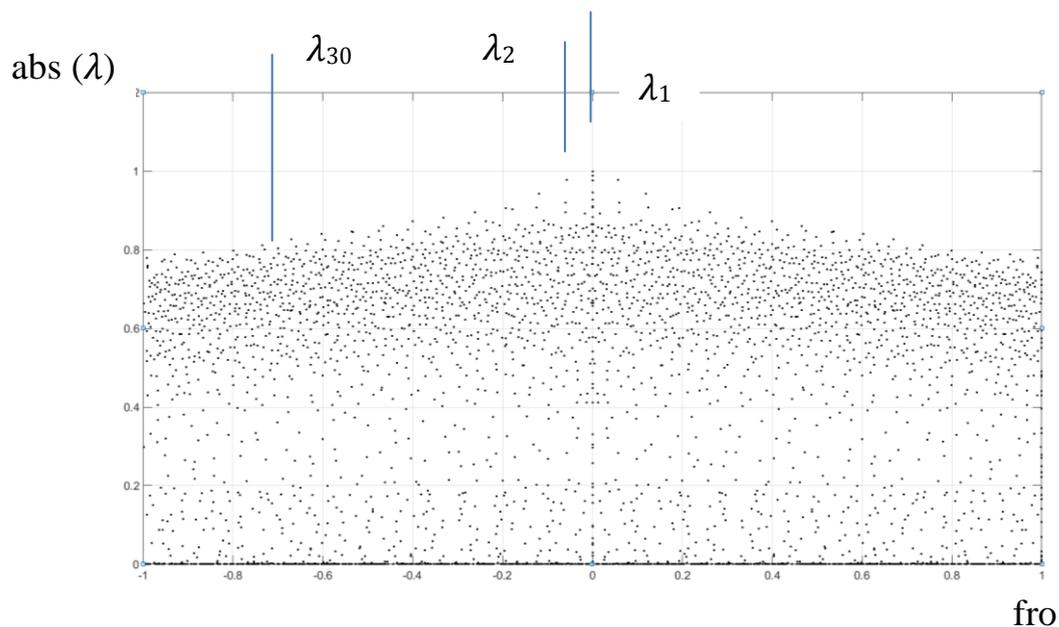

fro

Fig. 12 - Typical periods of fluctuations

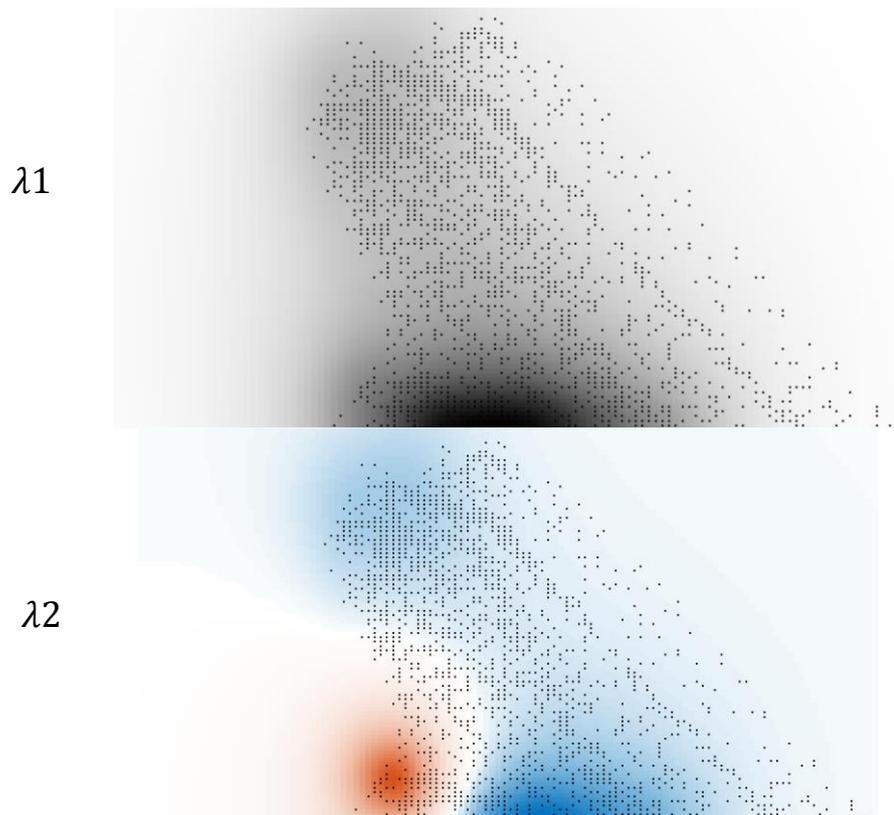
Fig. 13 – Human gait eigenforms in measured angles coordinates.

## Conclusions

We suggest the approach to construct and analyze a model of dynamic system by a sample signal, using Markov chains in the delay space. This class of models is in demand when it is impossible to describe a dynamic system by a system of differential equations.

For example, it seems a very complex task to build a control system for prostheses with a structural model of a control object, because we need to know moments of muscle forces, which are generated according to a complex law by the human nervous system.

On the other hand, the control system algorithm should be supplemented with a set of tools, helpful for the qualitative analysis of the resulting model. Among these tools, our article deals with the analysis of the dimensions of the resulting model, the attractor of the dynamic system and the characteristic frequencies of oscillations in the phase space.

# Bibliography

1. Coevoet E., Escande A., Duriez C. Soft robots locomotion and manipulation control using FEM simulation and quadratic programming // 2019 2nd IEEE International Conference on Soft Robotics (RoboSoft). - IEEE, 2019 .-- S. 739-745.
2. Volchek Yu. A. et al. Position of the artificial neural network model in medical expert systems // Juvenis scientia. - 2017. - No. 9.
3. Revuz D. Markov chains // M .: RFBR. - 1997.
4. Bogachev VI, Kirillov AI, Shaposhnikov SV On probabilistic and integrable solutions of the stationary Kolmogorov equation // Doklady Akademii Nauk. - Federal State Unitary Enterprise Academic Research and Publishing, Production, Printing and Book Distribution Center Nauka, 2011. - V. 438. - No. 2. - S. 154-159.
5. Nikulchev EV Geometric method for reconstructing systems from experimental data // Technical Physics Letters. - 2007. - T. 33. - No. 6. - S. 83-89.
6. Nefedova Yu.S., Shevtsova IG On non-uniform estimates of the rate of convergence in the central limit theorem // Theory of Probabilities and Its Applications. - 2012. - T. 57. - no. 1. - S. 62-97.
7. Chechulin V. L., Gracilev V. I. Stable regression estimation based on Chebyshev's inequality // ARTICLES IN THE JOURNAL "UNIVERSITY RESEARCH" 2009-2014. - 2015 . - S. 456-459.
8. Alpin Yu. A., Alpina VS The Perron – Frobenius theorem: proof using Markov chains. Zapiski nauchnykh seminars POMI. - 2008. - T. 359. - No. 0. - S. 5-16.
9. Anosov D. Differential equations: sometimes we solve, then we draw. - Litres, 2018.
10. Korotchenko MA, Burmistrov AV Modeling the dynamics of multiparticle ensembles using kinetic models // Educational resources and technologies. - 2016. - No. 2 (14).
11. Akimov OE Discrete mathematics: logic, groups, graphs // M .: Laboratory of basic knowledge. - 2001 .-- T. 352.
12. Svyatkina MN, Tarasov VB Third generation knowledge acquisition systems based on cognitive measurements // ALGORITHM FOR SEARCHING ANOMALIES IN PROCESSES BASED ON FUZZY TRENDS OF TIME SERIES. - 2014 .-- P. 58.
13. Staniswalis JG The kernel estimate of a regression function in likelihood-based models // Journal of the American Statistical Association. - 1989. - T. 84. - No. 405 .-- S. 276-283.

14. Dolgarev A.I., Some properties of the hypersphere of a space of dimension 6, in Theoretical and applied aspects of modern science. - 2014. - No. 5-1. - S. 23-27.
15. Nevsky M.V., On the geometric characteristics of an n-dimensional simplex, Modeling and Analysis of Information Systems, no. - 2011. - T. 18. - No. 2. - S. 52-64.
16. Sobol I. M. Method of Monte Carlo. - Science, 1985 .-- T. 46.